\documentclass[11pt]{article}
\usepackage{amsfonts}
\usepackage{amssymb}
\usepackage{amsmath}
\def\sqr#1#2{{\vcenter{\vbox{\hrule height.#2pt
              \hbox{\vrule width.#2pt height#1pt \kern#1pt \vrule
width.#2pt}
              \hrule height.#2pt}}}}
\def\signed #1{{\unskip\nobreak\hfil\penalty50
              \hskip2em\hbox{}\nobreak\hfil#1
              \parfillskip=0pt \finalhyphendemerits=0 \par}}
\def\endpf{\signed {$\sqr69$}}
\def\dbE{{\mathbb{E}}}
\def\dbF{{\mathbb{F}}}

\def\dbI{{\mathbb{I}}}

\def\dbL{{\mathbb{L}}}

\def\dbN{{\mathbb{N}}}

\def\dbP{{\mathbb{P}}}

\def\dbR{{\mathbb{R}}}

\def\a{\alpha}

\def\d{\delta}
\def\e{\varepsilon}
\def\z{\zeta}

\def\m{\mu}
\def\n{\nu}
\def\si{\sigma}

\def\f{\varphi}

\def\o{\omega}

\def\3n{\negthinspace \negthinspace \negthinspace }
\def\2n{\negthinspace \negthinspace }
\def\1n{\negthinspace }
\def\ns{\noalign{\smallskip} }
\def\ns{\noalign{\medskip} }
\def\ds{\displaystyle}
\def\D{\Delta}

\def\O{\Omega}
\def\cB{{\cal B}}

\def\cF{{\cal F}}

\def\cH{{\cal H}}

\def\cM{{\cal M}}

\def\cR{{\cal R}}
\def\cS{{\cal S}}

\def\cl{{\cal l}}
\def\no{\noindent}

\def\ss{\smallskip}
\def\ms{\medskip}
\def\bs{\bigskip}
\def\q{\quad}
\def\qq{\qquad}
\def\hb{\hbox}

\def\lan{\mathop{\langle}}
\def\ran{\mathop{\rangle}}

\def\esssup{\mathop{\rm esssup}}

\def\wt{\widetilde}
\def\cd{\cdot}
\def\cds{\cdots}

\def\as{\hbox{\rm a.s.{ }}}
\def\sgn{\hbox{\rm sgn$\,$}}

\def\cl{\overline}

\def\deq{\mathop{\buildrel\D\over=}}
\def\Re{{\mathop{\rm Re}\,}}
\def\Im{{\mathop{\rm Im}\,}}

\def\({\Big (}
\def\){\Big )}
\def\[{\Big[}
\def\]{\Big]}

\def\be{\begin{equation}}
\def\bel{\begin{equation}\label}
\def\ee{\end{equation}}
\def\bea{\begin{eqnarray}}
\def\eea{\end{eqnarray}}
\def\bt{\begin{theorem}}
\def\et{\end{theorem}}
\def\bc{\begin{corollary}}
\def\ec{\end{corollary}}
\def\bl{\begin{lemma}}
\def\el{\end{lemma}}
\def\bp{\begin{proposition}}
\def\ep{\end{proposition}}
\def\br{\begin{remark}}
\def\er{\end{remark}}
\def\ba{\begin{array}}
\def\ea{\end{array}}
\def\bd{\begin{definition}}
\def\ed{\end{definition}}

\newtheorem{lemma}{Lemma}[section]
\newtheorem{remark}{Remark}[section]

\newtheorem{theorem}{Theorem}[section]
\newtheorem{corollary}{Corollary}[section]

\newtheorem{definition}{Definition}[section]
\newtheorem{proposition}{Proposition}[section]

\makeatletter
   
   \@addtoreset{equation}{section}
\makeatother

\begin{document}

\title{\bf Representation of It\^o Integrals\\ by
Lebesgue/Bochner Integrals}
\author{Qi L\"{u}\thanks{School of Mathematics, Sichuan University, Chengdu 610064, China. {\small\it e-mail:} {\small\tt luqi59@163.com}. \ms},
  ~~~ Jiongmin Yong\thanks{Department of
Mathematics, University of Central Florida, FL 32816, USA. This work
was partially supported by the NSF under grant DMS-1007514.
{\small\it e-mail:} {\small\tt jyong@mail.ucf.edu}. \ms},~~~
  and~~~
Xu Zhang\thanks{Key Laboratory of Systems and Control, Academy of
Mathematics and Systems Sciences, Chinese Academy of Sciences,
Beijing 100190, China; and Yangtze Center of Mathematics, Sichuan
University, Chengdu 610064, China. This work is supported by the
NSFC under grants 10831007, 60821091 and 60974035,  and the
project MTM2008-03541 of the Spanish Ministry of Science and
Innovation. {\small\it e-mail:} {\small\tt xuzhang@amss.ac.cn}.}}
\maketitle
\begin{abstract}

In \cite{Yong 2004}, it was proved that as long as the integrand has
certain properties, the corresponding It\^o integral can be written
as a (parameterized) Lebesgue integral (or a Bochner integral). In
this paper, we show that such a question can be answered in a more
positive and refined way. To do this, we need to characterize the
dual of the Banach space of some vector-valued stochastic processes
having different integrability with respect to the time variable and
the probability measure. The later can be regarded as a variant of
the classical Riesz Representation Theorem, and therefore it will be
useful in studying other problems. Some remarkable consequences are
presented as well, including a reasonable definition of exact
controllability for stochastic differential equations and a
condition which implies a Black-Scholes market to be complete.

\end{abstract}

\bs

\no{\bf 2010 Mathematics Subject Classification}.  Primary 60G05;
Secondary 60H05, 60G07.

\bs

\no{\bf Key Words}. It\^o integral, Lesbegue integral, Bochner
integral, range inclusion, Riesz-type Representation Theorem.

\newpage

\section{Introduction}\label{s1}

Let $(\O,\cF,\dbF,\dbP)$ be a complete filtered probability space
with $\dbF=\{\cF_t\}_{t\ge0}$, on which a one-dimensional standard
Brownian motion $\{W(t)\}_{t\ge0}$ is defined so that $\dbF$ is its
natural filtration augmented by all the $\dbP$-null sets. Let $H$ be
a Banach space with the norm $|\cd|_H$ and with the dual space
$H^*$. For any $p\in [1,\infty)$, let $L_{\cF_T}^p(\O;H)$ be the set
of all $\cF_T$-measurable ($H$-valued) random variables $\xi:\O\to
H$ such that $\mathbb{E}|\xi|_{H}^p < \infty$. Next, for any
$p,q\in[1,\infty)$, put
\bel{1q1}\ba{ll}
\ns\ds L^p_\dbF(\O;L^q(0,T;H))=\Big\{\f:[0,T]\times\O\to
H\bigm|\f(\cd)\hb{
is $\dbF$-progressively measurable and}\\
\ns\ds\qq\qq\qq\qq\qq\qq\dbE\(\int_0^T|\f(t)|_H^qdt\)^{\frac{p}{q}}<\infty\Big\},\\
\ns\ds L^q_\dbF(0,T;L^p(\O;H))=\Big\{\f:[0,T]\times\O\to
H\bigm|\f(\cd)\hb{
is $\dbF$-progressively measurable and}\\
\ns\ds\qq\qq\qq\qq\qq\qq\int_0^T\(\dbE|\f(t)|_H^p\)^{\frac{q}
{p}}dt<\infty\Big\}.\ea\ee
In an obvious way, we may also define (for $1\le p,q<\infty$)
$$\left\{\ba{ll}
\ns\ds L^\infty_\dbF(\O;L^q(0,T;H)),\q
L^p_\dbF(\O;L^\infty(0,T;H)),\q L^\infty_\dbF
(\O;L^\infty(0,T;H)),\\
\ns\ds L^\infty_\dbF(0,T;L^p(\O;H)),\q
L^q_\dbF(0,T;L^\infty(\O;H)),\q
L^\infty_\dbF(0,T;L^\infty(\O;H)).\ea\right.$$
It is clear that
$$L^p_\dbF(\O;L^p(0,T;H))=L^p_\dbF(0,T;L^p(\O;H))\equiv
L^p_\dbF(0,T;H),\qq1\le p\le\infty.$$
Also, by Minkovski's inequality, it holds that
\bel{1.3}\left\{\ba{ll}
\ns\ds L^p_\dbF(\O;L^q(0,T;H))\subseteq
L^q_\dbF(0,T;L^p(\O;H)),\qq1\le p\le q\le\infty,\\
\ns\ds L^q_\dbF(0,T;L^p(\O;H))\subseteq
L^p_\dbF(\O;L^q(0,T;H)),\qq1\le q\le p\le\infty.\ea\right.\ee
In particular,
\bel{1.4}L^1_\dbF(0,T;L^p(\O;H))\subseteq
L^p_\dbF(\O;L^1(0,T;H)),\qq1\le p\le\infty.\ee

\ms

We now introduce two linear operators
\bel{I}\left\{\ba{ll}
\ns\ds\dbI:L^1_\dbF(\O;L^2(0,T;H))\to
L^1_{\cF_T}(\O;H)\q\hb{(when }H\hb{ is a Hilbert space)},\\
\ns\ds\dbI\big(\z(\cd)\big)=\int_0^T\z(t)dW(t),\qq\forall\;\z(\cd)\in
L^1_\dbF(\O;L^2(0,T;H)),\ea\right.\ee
and
\bel{L}\left\{\ba{ll}
\ns\ds\dbL:L^1_\dbF(0,T;H)\to L^1_{\cF_T}(\O;H),\\
\ns\ds\dbL\big(u(\cd)\big)=\int_0^Tu(t)dt,\qq\forall\; u(\cd)\in
L^1_\dbF(0,T;H).\ea\right.\ee
We call $\dbI$ and $\dbL$ the {\it It\^o integral operator} and the
{\it Lebesgue integral operator}, respectively. It is clear that
\bel{1.8}\left\{\ba{ll}
\ns\ds\dbI\(L^p_\dbF(\O;L^2(0,T;H))\)\subseteq
L^p_{\cF_T}(\O;H),\qq\forall\; p\in[1,\infty)\q\hb{(when }H\hb{ is a Hilbert space)},\\
\ns\ds\dbL\(L^p_\dbF(\O;L^1(0,T;H))\)\subseteq
L^p_{\cF_T}(\O;H),\qq\forall\; p\in[1,\infty).\ea\right.\ee

The first inclusion in (\ref{1.8}) can be refined (when $H$ is a
Hilbert space). Indeed, for any $\xi\in L^p_{\cF_T}(\O;H)$ (with
$p\in[1,\infty)$), $\dbE[\xi\,|\,\cF_t]$ is an $H$-valued continuous
$L^p$-martingale. Hence, by the Martingale Representation Theorem
(\cite{KS1988}), there is a unique $\z(\cd)\in
L^p_\dbF(\O;L^2(0,T;H))$ (called the {\it Malliavin derivative}
(\cite{Nualart1995}) of $\xi$ and sometimes denoted by $D_\cd\xi$)
such that
\bel{}\dbE[\xi\,|\,\cF_t]=\dbE\xi+\int_0^t\z(s)dW(s),\qq\forall\;
t\in[0,T].\ee
In particular, by taking $t=T$ in the above, we see that
\bel{}\xi=\dbE\xi+\int_0^T\z(s)dW(s).\ee
Therefore, in the case that $H$ is a Hilbert space, the first
inclusion in (\ref{1.8}) can be refined to the following
equality:
\bel{1.10}L^p_{\cF_T}(\O;H)=H\oplus\[\dbI\(L^p_\dbF(\O;L^2(0,T;H))\)\],\ee
where ``$\oplus$'' stands for a direct sum. Now, for the second
inclusion in (\ref{1.8}), we have the following simple result.

\ms

\bf Proposition 1.1. \sl Let $H$ be a Hilbert space and
$p\in[1,\infty)$. Then
\bel{12ooo3}\cl{\dbL\(L^p_\dbF(\O;L^1(0,T;H))\)}^{\,L^p_{\cF_T}(\O;H)}=L^p_{\cF_T}(\O;H),\ee
where $\cl G^{\,L^p_{\cF_T}(\O;H)}$ stands for the closure of $G$ in
$L^p_{\cF_T}(\O;H)$.

\ms

\it Proof. \rm For any $\z\in L^p_{\cF_T}(\O;H)$, let
$$\xi(t)=\dbE[\z\,|\,\cF_t],\qq t\in[0,T].$$
Then $\xi(\cd)$ is an $H$-valued $L^p$-martingale. By the martingale
representation theorem and the Burkholder--Davis--Gundy's
inequality, we have
$$\dbE\big|\xi(t)-\z|^p_H\le C\dbE\(\int_t^T|D_s\z|^2ds\)^{p\over2}\to0,\qq\hb{as }t\to T.$$
Now, for any $\d>0$, let
$$u_\d(t)={\xi(T-\d)\over\d}I_{[T-\d,T]}(t),\qq t\in[0,T].$$
Then $u_\d(\cd)\in L^p_\dbF(\O;L^\infty(0,T;H))\bigcap
L^\infty_\dbF(0,T;L^p(\O;H))\subseteq L^p_\dbF(\O;L^1(0,T;H))$, and
$$\dbE\Big|\int_0^Tu_\d(t)dt-\z\Big|^p=\dbE\big|\xi(T-\d)-\z\big|^p\to0,\qq\hb{as }
\d\to0,$$
proving the proposition. \endpf

\ms

\bf Remark 1.2. \rm From the proof of Proposition 1.1, it is easy to
see that we have proved the following stronger result than
(\ref{12ooo3}):
\bel{}
\cl{\dbL\(L^p_\dbF(\O;L^\infty(0,T;H))\)}^{\,L^p_{\cF_T}(\O;H)}=\cl{\dbL\(L^\infty_\dbF(0,T;L^p(\O;H))\)}^{\,L^p_{\cF_T}(\O;H)}=L^p_{\cF_T}(\O;H).
\ee

\ms

From Proposition 1.1, it is seen that we do not expect to have a
refinement for the Lebesgue integral operator $\dbL$ similar to
(\ref{1.10}). Instead, it is very natural for us to pose the
following problem:

\ms

{\bf Problem (E)} {\it Whether the following is true:
 \bel{1.11}
 \dbL\(L^p_\dbF(\O;L^1(0,T;H))\)=L^p_{\cF_T}(\O;H)\;?\ee}

Note that the above is equivalent to the following: When the range
of the operator $\dbL:L^p_\dbF(\O;L^1(0,T;H))\to L^p_{\cF_T}(\O;H)$
is closed. An interesting problem closely related to the above,
taking into account (\ref{1.10}), reads as follows.

\ms

{\bf Problem (R)} {\it Under what additional conditions on
$\z(\cd)\in L^p_\dbF(\O;L^2(0,T;H))$, there will be a $u(\cd)\in
L^p_\dbF(\O;L^1(0,T;H))$ such that the following holds
\bel{R}\int_0^T\z(t)dW(t)=\int_0^Tu(t)dt\;\qq\as?\ee}

\ss

For convenience, any $u(\cd)\in L^p_\dbF(\O;L^1(0,T;H))$ satisfying
(\ref{R}) is called a {\it representor} of $\z(\cd)$. Since the
It\^o integral in the usual sense can only be defined on Hilbert
spaces, we pose Problem (R) for the case that $H$ is a Hilbert
space. It is clear that when $u(\cd)$ is a representor of $\z(\cd)$
so is $u(\cd)+v(\cd)$ as long as $\ds\int_0^Tv(t)dt=0$, almost
surely. Therefore, if $\z(\cd)$ admits one representor, it admits
infinitely many representors. Problem (R) with $H=\dbR$ was posed
and studied in \cite{Yong 2004}. Various integrability conditions
were imposed on $\z(\cd)$ so that it admits a representor. Let us
now briefly recall several relevant results from \cite{Yong 2004},
which will give us some feelings about the representation (\ref{R}).
To this end, we define
\bel{u}u_\a(s)\equiv{1-\a\over{(T-s)^\a}}\int_0^s
{{\z(t)}\over{(T-t)^{1-\a}}}dW(t),\qq s\in[0,T),\ee
for $\a\in[0,1)$. The following is a summary of the relevant results
presented in \cite{Yong 2004}.

\ms

\bf Theorem 1.3. \sl {\rm(i)} Let $p\ge1$. For any $\z(\cd)\in
L^p_\dbF(\O;L^2(0,T;\dbR))$,
\bel{}u_0(\cd)\equiv\int_0^\cd{{\z(t)}\over{T-t}}dW(t)\in
\bigcup_{\e>0}L^p_\dbF(\O;L^2(0,T\2n-\e;\dbR)),\ee
and $(\ref{R})$ holds with $u(\cd)=u_0(\cd)$ in the following sense:
\bel{}\lim_{\e\to0}\dbE\left|\int_0^{T-\e}u_0(t)dt-\int_0^T\z(t)dW(t)\right|^p=0.\ee

{\rm(ii)} Suppose $\z(\cd)\in L^1_\dbF(0,T;L^2(\O;\dbR))$ such that
\bel{}\int_0^T\left[\int_0^s{{\dbE|\z(t)|^2}\over(T-t)^2}dt\right]^{1\over2}ds<\infty.\ee
Then
\bel{u0}u_0(\cd)\equiv\int_0^\cd{{\z(t)}\over{T-t}}dW(t)\in
L^1_\dbF(0,T;\dbR),\ee
and $(\ref{R})$ holds with $u(\cd)=u_0(\cd)$.

\ms

{\rm(iii)} Suppose $\z(\cd)\in L^1_\dbF(0,T;\dbR)$ such that for
some $\d>0$ the following holds:
\bel{}\int_0^T{{\dbE|\z(t)|^2}\over{(T-t)^\d}}dt<\infty.\ee
Then
\bel{}u_\a(\cd)\in
L^2_\dbF(\O;L^q(0,T;\dbR)),\qq\forall\;\a\in\hb{$({1-\d\over2},{1\over
q})$}\bigcap[0,1],\q q\in\hb{$[1,{2\over{2-\min(\d,1)}})$},\ee
and
\bel{}u_\a(\cd)\in
L^q_\dbF(0,T;\dbR),\qq\forall\;\a\in\hb{$(1-{\d\over2}-{1\over
q},{1\over q})$}\bigcap[0,1],\q
q\in\hb{$[1,{2\over{2-\min(\d,1)}})$}.\ee
Moreover, $(\ref{R})$ holds with $u(\cd)=u_\a(\cd)$.

\ms

{\rm(iv)} Suppose $\z(\cd)\in L^p_\dbF(0,T;\dbR)$ for some $p>2$.
Then
\bel{}u_\a(\cd)\in L^p_\dbF(\O;L^q(0,T;\dbR)),\qq\forall\;\a\in
\hb{$({1\over2},{1\over2}+{1\over p})$}\bigcap[0,1],\q
q\in\hb{$[1,{2p\over{p+2}})$}.\ee
Moreover, $(\ref{R})$ holds with $u(\cd)=u_\a(\cd)$.

\ms

\rm

The above shows that there are many $\z(\cd)\in
L^p_\dbF(\O;L^2(0,T;\dbR))$ such that one can find a corresponding
representor $u(\cd)$.

\ms

Note that although Problem (R) is posed for the case $H$ is a
Hilbert space, Problem (E) can be posed for general Banach space
since It\^o's integral is not involved here. The main purpose of
this paper is to give a positive answer to Problem (E) when $H$ is a
Banach space with $H^*$ having the Radon--Nikod\'{y}m property. Our
result seems to be a little surprising in some sense, and it refines
the results of \cite{Yong 2004} on Problem (R). More precisely, when
the answer to Problem (E) is positive, any $\z(\cd)\in
L^p_\dbF(\O;L^2(0,T;H))$ (when $H$ is a Hilbert space) admits a
representor $u(\cd)\in L^p_\dbF(\O;L^1(0,T;H))$, without assuming
further integrability conditions on $\z(\cd)$. This means that an
It\^o's integral on a given (fixed) interval can be represented by a
(parameterized) Bochner integral on that interval. We should
emphasize here that any representor $u(\cd)$ of $\z(\cd)\in
L^p_\dbF(\O;L^2(0,T;H))$ depends on $T$, in general. In another
word, it will be more proper to write
\bel{rep1}\int_0^T\z(t)dW(t)=\int_0^Tu(t,T)dt,\qq\as\ee
Hence, by allowing the upper limit to change, we should have
\bel{1.0124}\int_0^s\z(t)dW(t)=\int_0^s u(t,s)dt,\qq\forall\;
s\in[0,T],\q\as\ee
According to Theorem 1.3, when $\z(\cd)$ satisfies certain (better)
integrability conditions, we can find a representor of the following
form:
\bel{}u(t,s)={1-\a\over(s-t)^\a}\int_0^t{\z(r)\over(s-r)^{1-\a}}dW(r),\qq
0\le t<s\le T,\ee
for some $\a\in[0,1)$. Clearly, such an $s\mapsto u(t,s)$ is smooth
in $s\in(t,T]$. Therefore it is natural to further ask the following
question, without assuming the better integrability conditions on
$\z(\cd)$.

\ms

{\bf Problem (C)}  {\it For any $\z(\cd)\in
L^p_\dbF(\O;L^2(0,T;H))$, whether it has a representor $u(t,s)$
which is continuous with respect to the variable $s$?}

\ms

We will also show that the answer to Problem (C) is positive. Note
that, since the It\^o integral $\ds s\mapsto\int_0^s\z(t)dW(t)$ is
at most H\"older continuous up to order ${1\over2}$, generally, one
cannot expect that the differentiability of $s\mapsto u(t,s)$ (given
in (\ref{1.0124})). Nevertheless, it is natural to expect that
$s\mapsto u(t,s)$ is H\"older continuous up to order ${1\over2}$.
But, we do not have a proof for this yet.

\ms

\bf Remark 1.4. \rm The fact that $u(\cd)$ in (\ref{rep1})  depends
on $T$ tells us that, the positive answer to Problem (E) does not
mean that It\^o integrals can be completely replaced by
(parameterized) Bochner integrals.

\ms

The rest of this paper is organized as follows. In Section \ref{s2},
as a preliminary result, we establish a Riesz-type Representation
Theorem for the dual of the Banach space $L^p_\cM(X_1;L^q(X_2;H))$
(see Subsection \ref{sg1} for its definition). An interesting
byproduct in this section is the characterization on the dual of
$L^p_\dbF(\O;L^q(0,T;H))$ and $L^q_\dbF(0,T;L^p(\O;H))$, which will
be useful in some problems appeared in stochastic distributed
parameter control systems and/or stochastic partial differential
equations. Section \ref{s3} is addressed to giving answers to
Problems (E) and (R). Section \ref{s4} is devoted to answering
Problem (C), for which the key tool we employ is the continuous
selection theorem in \cite{Micheal}. In Section \ref{s5}, we present
two remarkable consequences of our positive solution to Problem (E),
one of which is related to the reasonable formulation of exact
controllability for stochastic differential equations, and the other
a condition to guarantee a Black-Scholes market to be complete.

\section{The Dual of $L^p_\cM(X_1;L^q(X_2;H))$}\label{s2}

As a key preliminary to answer Problem (E), we need to characterize
the dual of $L^p_\dbF(\O;L^q(0,T;H))$ and $L^q_\dbF(0,T;L^p(\O;H))$.
We will go a little further by considering the dual of
$L^p_\cM(X_1;L^q(X_2;H))$, which will be be defined below. It seems
to us that this result has its own interest.

\subsection{Statement of the result}\label{sg1}

Let $(X_1,\cM_1,\m_1)$ and $(X_2,\cM_2,\m_2)$ be two finite
measure spaces. Let $\cM$ be a sub-$\si$-field of
$\cM_1\otimes\cM_2$ (the $\si$-field generated by
$\cM_1\times\cM_2$), and for any $1\le p,q<\infty$, let
$$\ba{ll}
\ns\ds L^p_\cM(X_1;L^q(X_2;H))=\Big\{\f:X_1\times X_2\to
H\bigm|\f(\cd)\hb{ is
$\cM$-measurable and}\\
\ns\ds\qq\qq\qq\qq\qq\qq\qq\qq\int_{X_1}\left(\int_{X_2}|\f(x_1,x_2)|_H^qd\m_2\right)^{p\over
q}d\m_1<\infty\Big\}.\ea$$
Likewise, let
$$\ba{ll}
\ns\ds L^\infty_\cM(X_1;L^q(X_2;H))=\Big\{\f:X_1\times X_2\to
H\bigm|\f(\cd)\hb{ is
$\cM$-measurable and}\\
\ns\ds\qq\qq\qq\qq\qq\qq\qq\qq\esssup_{x_1\in
X_1}\left(\int_{X_2}|\f(x_1,x_2)|_H^qd\m_2\right)^{1\over
q}<\infty\Big\},\ea$$
$$\ba{ll}
\ns\ds L^p_\cM(X_1;L^\infty(X_2;H))=\Big\{\f:X_1\times X_2\to
H\bigm|\f(\cd)\hb{ is
$\cM$-measurable and}\\
\ns\ds\qq\qq\qq\qq\qq\qq\qq\qq\int_{X_1}\left(\esssup_{x_2\in
X_2}|\f(x_1,x_2)|_H^p\right)d\m_1<\infty\Big\},\ea$$
$$\ba{ll}
\ns\ds L^\infty_\cM(X_1;L^\infty(X_2;H))=\Big\{\f:X_1\times X_2\to
H\bigm|\f(\cd)\hb{ is
$\cM$-measurable and}\\
\ns\ds\qq\qq\qq\qq\qq\qq\qq\qq\esssup_{(x_1,x_2)\in X_1\times
X_2}|\f(x_1,x_2)|_H<\infty\Big\}.\ea$$
We denote
$$L^p_\cM(X_1\times X_2;H)=L^p_\cM(X_1;L^p(X_2;H)),\qq1\le p\le\infty.$$
Also, for any $f\in L^p_\cM(X_1;L^q(X_2;H))$ ($1\le
p,q\le\infty$), we denote
\bel{norm}\|f\|_{p,q,H}\equiv\|f\|_{L^p_\cM(X_1;L^q(X_2;H))}\deq\left[\int_{X_1}\left(\int_{X_2}
|f(x_1,x_2)|_H^qd\m_2\right)^{p\over q}d\m_1\right]^{1\over p}.\ee
The definition of $\|f\|_{\infty,q,H}$, $\|f\|_{p,\infty,H}$ and
$\|f\|_{\infty,\infty,H}$ are obvious. Let
\bel{}\|f\|_{p,H}\equiv\|f\|_{p,p,H}\;,\qq\qq1\le p\le\infty.\ee
The definition of $L^q_\cM(X_2;L^p(X_1;H))$ ($1\le p,q\le\infty$) is
similar. By H\"older's inequality and Minkovski's inequality, we
have the following inclusions:
\bel{basd}\ba{ll}
\ns\ds L^p_\cM(X_1;L^q(X_2;H))\subseteq
L^r_\cM(X_1;L^s(X_2;H)),\qq1\le r\le p\le\infty,\q1\le s\le
q\le\infty,\ea\ee
and (comparing with (\ref{1.3})--(\ref{1.4})),
\bel{}\left\{\ba{ll}
\ns\ds L^p_\cM(X_1;L^q(X_2;H))\subseteq
L^q_\cM(X_2;L^p(X_1;H)),\qq1\le
p\le q\le\infty,\\
\ns\ds L^p_\cM(X_1;L^q(X_2;H))\supseteq
L^q_\cM(X_2;L^p(X_1;H)),\qq1\le q\le p\le\infty.\ea\right.\ee
Next, for any $p\in[1,\infty]$, denote
$$p'=\left\{\ba{ll}
\ns\ds{p\over p-1},\qq1<p<\infty,\\
\ns\ds1,\qq\qq p=\infty,\\
\ns\ds\infty,\qq\q~ p=1.\ea\right.$$
The definition of $q'\in[1,\infty]$ for $q\in[1,\infty]$ is similar.
We have the following result.

\ms

\bf Lemma 2.1. \sl Let $H$ be a Banach space, $(X_1,\cM_1,\m_1)$ and
$(X_2,\cM_2,\m_2)$ be two finite measure spaces, $\cM$ be a
sub-$\si$-field of $\cM_1\otimes\cM_2$, and let $1\le p,q<\infty$.
Then,  $H^*$ has the Radon--Nikod\'ym property with respect to
$(X_1\times X_2, \cM, \mu_1\times\mu_2)$ if and only if for any
$F\in L^p_\cM(X_1;L^q(X_2;H))^*$, there exists a unique $g\in
L^{p'}_\cM(X_1;L^{q'}(X_2;H^*))$ such that
\bel{}F(f)=\int_{X_1\times
X_2}(f(x_1,x_2),g(x_1,x_2))_{H,H^*}d\m_1d\m_2,\qq\forall\; f\in
L^p_\cM(X_1;L^q(X_2;H)),\ee
and
\bel{F=g}\|F\|_{L^p_\cM(X_1;L^q(X_2;H))^*}=\|g\|_{p',q',H^*}.\ee

\ms

\rm

Due to the above result, we make the following identification (for
the case that $H^*$ has the Radon--Nikod\'ym property with respect
to $(X_1\times X_2, \cM, \mu_1\times\mu_2)$):
\bel{}L^p_\cM(X_1;L^q(X_2;H))^*=L^{p'}_\cM(X_1;L^{q'}(X_2;H^*)),\qq1\le
p,q<\infty.\ee
The above is a Riesz-type Representation Theorem for the dual of
space $L^p_\cM(X_1;L^q(X_2;H))$. It seems to us that Lemma 2.1
should be a known result but we have not found an exact reference.
Therefore, for reader's convenience, we provide a detailed proof in
the next three subsections. As a corollary of Lemma 2.1, we will
characterize the dual of $L^p_\dbF(\O;L^q(0,T;H))$ and
$L^q_\dbF(0,T;L^p(\O;H))$ in the last subsection.

\ms

The main idea for the proof of Lemma 2.1 is similar to that of the
relevant result in \cite[Appendix B, pp. 375--376]{C} (see also
\cite[Theorem 1, Chapter IV, pp. 98--99]{DU}). However, Lemma 2.1
does not follow from the main result in \cite[Appendix B]{C} because
the later considered only the special case that $p=q$ and $H=\dbR$,
for which, by Fubini's Theorem, one can reduce the problem to the
case with one measure on the product space. Also, Lemma 2.1 does not
seem to be a corollary of \cite[Theorem 1, Chapter IV, pp.
98--99]{DU} because of the very fact that our $\cM$ is an
``interconnecting" sub-$\si$-field of the $\si$-field generated by
$\cM_1\times\cM_2$.

\subsection{Proof of the necessity in Lemma 2.1 for the case
$H=\dbR$}\label{2s2}

As a key step to prove Lemma 2.1, in this subsection we show first
the ``only if" part of this lemma for the special case $H=\dbR$.

For any $g\in L^{p'}_\cM(X_1;L^{q'}(X_2;\dbR))$, define $F_g:
L^p_\cM(X_1;L^q(X_2;\dbR))\mapsto\dbR$ by
$$F_g(f)=\int_{X_1\times X_2}f(x_1,x_2)g(x_1,x_2)d\m_1d\m_2,\qq\forall\;f\in
L^p_\cM(X_1;L^q(X_2;\dbR)).$$
By the linearity of the integral, $g\mapsto F_g$ is a linear map. It
follows from H\"older's inequality that
$$|F_g(f)|\le\|f\|_{p,q,\dbR}\|g\|_{p',q',\dbR},\qq\forall\;
f\in L^p_\cM(X_1;L^q(X_2;\dbR)).$$
Hence $F_g\in L^p_\cM(X_1;L^q(X_2;\dbR))^*$ and
\bel{F<g}\|F_g\|_{L^p_\cM(X_1;L^q(X_2;\dbR))^*}\le\|g\|_{p',q',\dbR}.\ee
Therefore, $g\mapsto F_g$ is a linear non-expanding map. Now, we
show that this map is surjective and is an isometry.

To show the surjectivity of $g\mapsto F_g$, take any $F\in
L^p_\cM(X_1;L^q(X_2;\dbR))^*$. Since for any $A\in\cM$, $I_A\in
L^p_\cM(X_1;L^q(X_2;\dbR))$, we may define
$$\nu(A)=F(I_A),\qq\forall\;A\in\cM.$$
Then $\nu$ is a totally finite signed measure on $(X_1\times
X_2,\cM)$, and $\n<\3n<\m_1\times\m_2$. By the Radon-Nikod\'ym
Theorem, there is an $\cM$-measurable map $g\in L^1_\cM(X_1\times
X_2;\dbR)$ such that
$$\nu(A)=\int_A gd\m_1d\m_2,\qq\forall\; A\in\cM,$$
i.e.,
$$F(I_A)=\int_{X_1\times X_2}gI_Ad\m_1d\m_2,\qq\forall\; A\in
\cM.$$
Consequently, for any $\cM$-measurable simple functions $f$,
$$F(f)=\int_{X_1\times X_2}f(x_1,x_2)g(x_1,x_2)d\m_1d\m_2.$$
Select a sequence $\{A_n\}_{n=1}^\infty\subset\cM$ such that
\bel{}A_n\subset A_{n+1},\q
n=1,2,\cds,\qq(\m_1\times\m_2)\left(\big(X_1\times X_2\big)\setminus
\bigcup_{n=1}^\infty A_n\right)=0,\ee
and $g$ is bounded on each $A_n$. For any $n\ge1$, note that
$$f\mapsto\int_{X_1\times
X_2}f(x_1,x_2)g(x_1,x_2)I_{A_n}(x_1,x_2)d\m_1d\m_2$$
is a bounded linear functional on $L^p_\cM(X_1;L^q(X_2;\dbR))$ which
agrees with $F$ on all $\cM$-measurable simple functions which
vanishes off $A_n$. It follows that
\bel{5.3}\ba{ll}
\ns\ds F(fI_{A_n})=\int_{X_1\times
X_2}fgI_{A_n}d\m_1d\m_2,\qq\forall\; f\in
L^p_\cM(X_1;L^q(X_2;\dbR)).\ea\ee
Since $gI_{A_n}$ is bounded, one has $gI_{A_n}\in
L^{p'}_\cM(X_1;L^{q'}(X_2;\dbR))$. We claim that $g\in
L^{p'}_\cM(X_1;L^{q'}(X_2;\dbR))$, and
\bel{F>g}\|g\|_{p',q';\dbR}\le\|F\|_{L^p_\cM(X_1;L^q(X_2;\dbR))^*}.\ee
To show this, we distinguish four cases.

\ms

{\bf Case 1: $p,q\in(1,\infty)$}. Choose
$$f=\left\{\ba{ll}
\ns\ds a\left(\int_{X_2}|g|^{q'}I_{A_n}d\m_2\right)^{{p'\over
q'}-1}|g|^{q'-1}(\sgn g)I_{A_n},\q &\hb{ if }\ds\int_{X_2}|g|^{q'}
I_{A_n}d\m_2\ne0,\\
\ns\ds0,&\hb{ if }\ds\int_{X_2}|g|^{q'}I_{A_n}d\m_2=0,\ea\right.$$
where
$$a=\left[\int_{X_1}\left(\int_{X_2}|g|^{q'}I_{A_n}d\m_2\right)^{p'\over q'}d\m_1\right]^{{1\over p'}-1}.$$
Then
$$\ba{ll}
\ns\ds\|f\|_{p,q}=\left[\int_{X_1}\left(\int_{X_2}|f|^qd\m_2\right)^{p\over
q}d\m_1\right]^{1\over p}\\
\ns\ds=\left\{\int_{X_1}\left[\int_{X_2}a^q\left(\int_{X_2}|g|^{q'}I_{A_n}d\m_2\right)^{\({p'\over
q'}-1\)q}|g|^{(q'-1)q}I_{A_n}d\m_2\right]^{p\over
q}d\m_1\right\}^{1\over p}\\
\ns\ds=a\left\{\int_{X_1}\left[\left(\int_{X_2}|g|^{q'}I_{A_n}d\m_2\right)^{\({p'\over
q'}-1\)p}\left(\int_{X_2}|g|^{q'}I_{A_n}d\m_2\right)^{p\over
q}\right]d\m_1\right\}^{1\over p}\\
\ns\ds=a\left\{\int_{X_1}\left(\int_{X_2}|g|^{q'}I_{A_n}d\m_2\right)^{p'\over
q'}d\m_1\right\}^{1\over p}=1.\ea$$
Taking the above $f$ in (\ref{5.3}), we find that
$$\ba{ll}
\ns\ds F(f)=\int_{X_1}\int_{X_2}fgI_{A_n}d\m_2d\m_1
=a\int_{X_1}\left[\int_{X_2}\left(\int_{X_2}|g|^{q'}I_{A_n}d\m_2\right)^{{p'
\over
q'}-1}|g|^{q'}I_{A_n}d\m_2\right]d\m_1\\
\ns\ds\qq\qq=a\int_{X_1}\left(\int_{X_2}|g|^{q'}I_{A_n}d\m_2\right)^{p'\over
q'}d\m_1=\left[\int_{X_1}\left(\int_{X_2}|g|^{q'}I_{A_n}d\m_2\right)^{p'\over
q'}d\m_1\right]^{1\over p'}\\
\ns\ds\qq\qq=\|gI_{A_n}\|_{p',q';\dbR},\ea$$
which gives
$$\|gI_{A_n}\|_{p',q';\dbR}\le\|F\|_{L^p_\cM(X_1;L^q(X_2;\dbR))^*}.$$
Letting $n\to\infty$, by making use of Fatou's Lemma, one concludes
(\ref{F>g}).

\ms

{\bf Case 2: $p=1$, $1<q<\infty$}. In this case, we first take
$p\in(1,\infty)$, and take $f$ as in Case 1. Then
$$\|f\|_{1,q}=\int_{X_1}\left(\int_{X_2}|f|^qd\m_2\right)^{1\over
q}d\m_1\le\[\int_{X_1}\left(\int_{X_2}|f|^qd\m_2\right)^{p\over
q}d\m_1\]^{1\over p}\m_1(X_1)^{1\over p'}=\m_1(X_1)^{1\over p'}.$$
Consequently,
$$\|gI_{A_n}\|_{p',q';\dbR}=F(f)\le\|F\|_{L^1_\cM(X_1;L^q(X_2;\dbR))^*}
\|f\|_{1,q}\le\|F\|_{L^1_\cM(X_1;L^q(X_2;\dbR))^*}\m_1(X_1)^{1\over
p'}.$$
Letting $n\to\infty$ and then letting $p\to1$ (which means
$p'\to\infty$), we obtain
\bel{}\|g\|_{\infty,q';\dbR}\le\|F\|_{L^1_\cM(X_1;L^q(X_2;\dbR))^*},\ee
which is (\ref{F>g}) for the case $p=1$.

\ms

{\bf Case 3: $1<p<\infty$, $q=1$}. In this case, we first take
$q\in(1,\infty)$, and take $f$ as in Case 1. Then
$$\|f\|_{p,1;\dbR}=\left[\int_{X_1}\left(\int_{X_2}|f|d\m_2\right)^pd\m_1\right]^{1\over
p}\le\left[\int_{X_1}\left(\int_{X_2}|f|^qd\m_2\right)^{p\over
q}d\m_1\right]^{1\over p}\m_2(X_2)^{1\over q'}=\m_2(X_2)^{1\over
q'}.$$
Hence,
$$\ba{ll}
\ns\ds\|gI_{A_n}\|_{p',q';\dbR}=F(f)\le\|F\|_{L^p_\cM(X_1;L^1(X_2;\dbR))^*}
\|f\|_{p,1;\dbR}\le\|F\|_{L^p_\cM(X_1;L^1(X_2;\dbR))^*}\m_2(X_1)^{1\over
q'}\ea$$
Letting $n\to\infty$ and then letting $q\to1$ (which means
$q'\to\infty$), we obtain
\bel{}\|g\|_{p',\infty;\dbR}\le\|F\|_{L^p_\cM(X_1;L^1(X_2;\dbR))^*},\ee
which is the case of (\ref{F>g}) for $q=1$.

\ms

{\bf Case 4: $p=q=1$}. In this case, we still first let
$p,q\in(1,\infty)$, and take $f$ as in Case 1 with $q=r$. Then
$$\ba{ll}\ds\|f\|_{1,1}=\int_{X_1}\int_{X_2}|f|d\m_2d\m_1\\
\ns\ds\le\left[\int_{X_1}\left(\int_{X_2}|f|^qd\m_2\right)^{p\over
q}d\m_1\right]^{1\over p}\m_1(X_1)^{1\over p'}\m_2(X_2)^{1\over
q'}=\m_1(X_1)^{1\over p'}\m_2(X_2)^{1\over q'}.\ea$$
Consequently,
$$\ba{ll}
\ns\ds\|gI_{A_n}\|_{p',q';\dbR}=F(f)\le\|F\|_{L^1_\cM(X_1;L^1(X_2;\dbR))^*}
\|f\|_{1,1}\le\|F\|_{L^1_\cM(X_1;L^1(X_2;\dbR))^*}\m_1(X_1)^{1\over
p'}\m_2(X_1)^{1\over q'}.\ea$$
Letting $n\to\infty$ and then letting $p,q\to1$ (which means
$p',q'\to\infty$), we obtain
\bel{}\|g\|_{\infty;\dbR}\le\|F\|_{L^1_\cM(X_1;L^1(X_2;\dbR))^*},\ee
which is the case of (\ref{F>g}) for $p,q=1$.

\ms

Finally, (\ref{F<g}) means that
$F_g\in(L^p_\cM(X_1;L^q(X_2;\dbR)))^*$ and since $F$ and $F_g$
coincides on a dense subset of $L^p_\cM(X_1;L^q(X_2;\dbR))$, one has
$F=F_g$. Also, (\ref{F=g}) follows easily from (\ref{F<g}) and
(\ref{F>g}).

\subsection{Proof of the necessity in Lemma 2.1 for the general case}

We are now in a position to prove the ``only if" part of Lemma 2.1
for the general case. The proof is divided into two steps.

\ms

{\it Step 1}. We show that $ L_{\cM}^{p'}(X_1;L^{q'}(X_2;H^*))$ is
isometrically isomorphic to a subspace $\cal H$ of
$L_{\cM}^p(X_1;L^q(X_2;H))^*$.

For any given $g\in L_{\cM}^{p'}(X_1;L^{q'}(X_2;H^*))$, define a
linear functional $F_g$ on $L_{\cM}^p(X_1;L^q(X_2;H))$ as follows:
\begin{equation}
F_g(f) = \int_{X_1\times X_2}\langle f(x_1,x_2),g(x_1,x_2)
\rangle_{H,H^*}d\mu_1 d\mu_2,\qq \forall f \in
L_{\cM}^p(X_1;L^q(X_2;H)).
\end{equation}
Then, by means of the H\"older inequality and similar to
(\ref{F<g}), we conclude that $F_g$ belongs to
$L_{\cM}^p(X_1;L^q(X_2;H))^*$, and
\begin{equation}\label{boundFg}
\|F_g\|_{L_{\cM}^p(X_1;L^q(X_2;H))^*} \leq \|g\|_{p'q',H^*}.
\end{equation}
Therefore the norm of $F_g$ is not greater than $\|g\|_{p'q',H^*}$.
Define
$$\cH\equiv \{F_g\;|\;g\in L_{\cM}^{p'}(X_1;L^{q'}(X_2;H^*))\}.$$

It remains to prove the reverse of inequality (\ref{boundFg}).
Clearly, without loss of generality, we may assume that $g\not=0$.

Suppose first that $\ds g=\sum_{i=1}^{\infty}h_i^*I_{E_i}$ where
$h_i^*$ is a sequence in $H^*$ and $\{E_i\}_{i=1}^{\infty}$ is a
countable partition of $X_1\times X_2$ by members of $\cM$ with
$(\mu_1\times\mu_2)(E_i)>0$ for all $i$. Since we have shown that
$L^p_\cM(X_1;L^q(X_2;\dbR))^*=L^{p'}_\cM(X_1;L^{q'}(X_2;\dbR))$ (in
Subsection \ref{2s2}) and noting that $0<|g|_{H^*}\in
L^{p'}_\cM(X_1;L^{q'}(X_2;\dbR))$, for any $\e>0$, there exists a
nonnegative function $\varphi\in L_{\cM}^p(X_1;L^q(X_2;\mathbb{R}))$
such that
$$0<\|\varphi\|_{p,q}\le1,\qq\|g\|_{p'q',H^*}-\e
\leq \int_{X_1\times X_2}|g|_{H^*}\varphi d\mu_1 d\mu_2.$$
Further, choose $h_i\in H$ with $|h_i|_{H}=1$ such that
$$|h_i^*|_{H^*} -\frac{\e}{\|\varphi\|_{1,1}}\leq h_i^*(h_i),$$
and define
$$f=\sum_{i=1}^\infty \varphi h_iI_{E_i}\in L_{\cM}^p(X_1;L^q(X_2;H)).$$
Then we have that $\|f\|_{p,q,H} = \|\varphi\|_{p,q} \leq 1$, and we
have that
$$\ba{ll}
\ns\ds\int_{X_1\times X_2}\lan f(x_1,x_2),g(x_1,x_2)
\rangle_{H,H^*}d\mu_1 d\mu_2= \int_{X_1\times X_2} \varphi
\sum_{i=1}^\infty\lan h_i,h_i^*
\rangle_{H,H^*}\chi_{E_i}d\mu_1 d\mu_2 \\
\ns\ds\ge\int_{X_1\times X_2}\varphi\sum_{i=1}^\infty \Big(
|h_i^*|_{H^*}-\frac{\e}{\|\varphi\|_{1,1}} \Big)\chi_{E_i}d\mu_1 d\mu_2 \\
\ns\ds\ge \int_{X_1\times X_2}|g|_{H^*}\varphi d\mu_1 d\mu_2 -
\frac{\e}{\|\varphi\|_{1,1}}\int_{X_1\times X_2}\varphi d\mu_1
d\mu_2\ge\|g\|_{p',q',H^*} - 2\e.\ea$$
This gives
$$\|F_g\|_{L_{\cM}^p(X_1;L^q(X_2;H))^*} \geq \|g\|_{p'q',H^*},$$
and therefore
$$\|F_g\|_{L_{\cM}^p(X_1;L^q(X_2;H))^*} =
\|g\|_{p'q',H^*},$$
whenever $g\in L_{\cM}^{p'}(X_1;L^{q'}(X_2;H^*))$ is countably
valued.

For the general case, we choose a sequence
$\{g_n\}_{n=1}^\infty\subset L_{\cM}^{p'}(X_1;L^{q'}(X_2;H^*))$ such
that each $g_n$ is countably valued  and
\begin{equation}\label{gng}
\lim_{n\to\infty}\|g_n - g\|_{p',q',H^*} = 0.
\end{equation}
We have obtained that
\begin{equation}\|F_{g_n}\|_{L_{\cM}^p(X_1;L^q(X_2;H))^*}
=\|g_n\|_{p',q',H^*},\nonumber
\end{equation}
and by virtue of (\ref{boundFg}),
\begin{equation}
\|F_{g_n}-F_{g}\|_{L_{\cM}^p(X_1;L^q(X_2;H))^*}=\|F_{g_n-g}\|_{L_{\cM}^p(X_1;L^q(X_2;H))^*}
\le\|g_n-g\|_{p',q',H^*}.\nonumber
\end{equation}
Therefore, noting (\ref{gng}), we end up with
\begin{eqnarray*}
\|F_g\|_{L_{\cM}^p(X_1;L^q(X_2;H))^*}
=\lim_{n\to\infty}\|F_{g_n}\|_{L_{\cM}^p(X_1;L^q(X_2;H))^*} =
\lim_{n\to\infty}\|g_n\|_{p'q',H^*}=\|g\|_{p'q',H^*}.
\end{eqnarray*}
Hence we get that $ L_{\cM}^{p'}(X_1;L^{q'}(X_2;H^*))$ is
isometrically isomorphic to $\cal H$.

\ms

{\it Step 2}. We show that the subspace $\cal H$ is equal to
$L_{\cM}^p(X_1;L^q(X_2;H))^*$.

\ms

To this end, for $F\in L_{\cM}^p(X_1;L^q(X_2;H))^*$, we define
\begin{equation}\label{Gdsqq}
G(E)(h) = F(hI_{E}),\qq \forall\; E\in \cM, \ h\in H.
\end{equation}
By
$$|F(hI_{E})|\le\|F\|_{L_{\cM}^p(X_1;L^q(X_2;H))^*}\|hI_{E}\|_{p,q,H}\le
\|F\|_{L_{\cM}^p(X_1;L^q(X_2;H))^*}|h|_{H}\|I_{E}\|_{p,q},$$
we see that $G:\cM\to H^*$ and it is countably additive. Let
$E_1,\cdots,E_n$ ($n\in\dbN$) be a partition of $X_1\times X_2$ by
members of $\cM$ with $(\mu_1\times\mu_2)(E_i)>0$ for all $1\leq i
\leq n$. Then $G(E_i)\in H^*$. Define
$$G_{E_i}^1(h)=\Re G(E)(h),\q G_{E_i}^2(h)=\Im G(E)(h),\qq\forall
h\in H.$$
Clearly, $|G(E_i)|_{H^*}\le
|G_{E_i}^1|_{H^*}+|G_{E_i}^2|_{H^*}$. Noting that both $G_{E_i}^1$
and $G_{E_i}^2$ are real functionals, we see that, for any $\e>0$,
one can find $h_i^1$ and $h_i^2$ in the closed unit ball of $H$ such
that
 $$
 |G_{E_i}^1|_{H^*}-\frac{\e}{2n}<\Re G(E_i)(h_i^1),\qq |G_{E_i}^2|_{H^*}-\frac{\e}{2n}<\Im G(E_i)(h_i^2).
 $$
It follows that
\begin{eqnarray*}
&\,& \sum_{i=1}^n|G(E_i)|_{H^*}-\e<\Re\sum_{i=1}^n G(E_i)(h_i^1)+\Im\sum_{i=1}^n G(E_i)(h_i^2) \\
&\,& \qq\qq\qq\qq= \Re F\(\sum_{i=1}^n h_i^1I_{E_i}\)+\Im F\(\sum_{i=1}^n h_i^2I_{E_i}\)\\
&\,& \qq\qq\qq\qq \leq \|F\|_{L_{\cM}^p(X_1;L^q(X_2;H))^*}
\left(\Big\|\sum_{i=1}^n h_i^1I_{E_i}\Big\|_{p,q,H}
+\Big\|\sum_{i=1}^n h_i^2I_{E_i}\Big\|_{p,q,H} \right)\\
&\,& \qq\qq\qq\qq \leq
2\|F\|_{L_{\cM}^p(X_1;L^q(X_2;H))^*}\Big\|\sum_{i=1}^n
I_{E_i}\Big\|_{p,q}\\
&\,& \qq\qq\qq\qq \leq 2\|F\|_{L_{\cM}^p(X_1;L^q(X_2;H))^*}
\mu_1(X_1)^{\frac{1}{p}}\mu_2(X_2)^{\frac{1}{q}}.
\end{eqnarray*}
Hence $|G(X_1\times X_2)|_{H^*}< \infty$ and $G$ is a
$(\mu_1\times\mu_2)$-continuous vector-valued measure of bounded
variation. Since $H^*$ has the Radon-Nikod\'{y}m property with
respect to $(X_1\times X_2, \cM, \mu_1\times\mu_2)$, there exists a
Bochner integrable $g:X_1\times X_2 \to H^*$ such that
\begin{equation}\label{GE}
G(E) = \int_E g d\mu_1 d\mu_2,\qq  \forall\; E\in \cM.
\end{equation}

Clearly, if $f\in L_{\cM}^p(X_1;L^q(X_2;H))$ is a simple function,
then
$$
F(f) = \int_{X_1\times X_2}\langle f(x_1,x_2),g(x_1,x_2)
\rangle_{H,H^*}d\mu_1 d\mu_2.
$$
Select an expanding sequence $\{E_n\}_{n=1}^\infty$ in $\cM$ such
that $\ds\bigcup_{n=1}^\infty E_n = X_1\times X_2$ and such that $g$
is bounded on each $E_n$. Fixing arbitrarily an $n_0\in\dbN$ and
noting that $\ds\int_{E_{n_0}}\langle \cdot,g(x_1,x_2)
\rangle_{H,H^*}d\mu_1 d\mu_2$ is a bounded linear functional on
$L_{\cM}^p(X_1;L^q(X_2;H))$ which agrees with $F$ on all simple
functions supported on $E_{n_0}$, it follows that
\begin{equation}\label{opqq}
F(fI_{E_{n_0}}) = \int_{X_1\times X_2}\langle
f(x_1,x_2),g(x_1,x_2)I_{E_{n_0}} \rangle_{H,H^*}d\mu_1 d\mu_2,\q
\forall\; f\in L_{\cM}^p(X_1;L^q(X_2;H)).
\end{equation}
Further, since $gI_{E_{n_0}}$ is bounded, one has $gI_{E_{n_0}}\in
L_{\cM}^{p'}(X_1;L^{q'}(X_2;H^*))$ and
\begin{equation}\label{esti1}\|gI_{E_{n_0}}\|_{p',q',H^*} \leq
\|F\|_{L_{\cM}^p(X_1;L^q(X_2;H))^*}.
\end{equation}
Since inequality (\ref{esti1}) holds for each $n_0$, by the Monotone
Convergence Theorem, we conclude that $g\in
L_{\cM}^{p'}(X_1;L^{q'}(X_2;H^*))$.

Finally, for any $f\in L_{\cM}^p(X_1;L^q(X_2;H))$, it follows from
(\ref{opqq}) that
$$\ba{ll}
\ns\ds F(f)=\lim_{n\to\infty} \int_{X_1\times X_2}\langle
f(x_1,x_2),g(x_1,x_2)I_{E_{n}} \rangle_{H,H^*}d\mu_1
d\mu_2\\
\ns\ds\qq\;=\int_{X_1\times X_2}\langle f(x_1,x_2),g(x_1,x_2)
\rangle_{H,H^*}d\mu_1 d\mu_2=F_g(f).\ea$$
This means that $F=F_g$. Hence $ L_{\cM}^p(X_1;L^q(X_2;H))^*$
coincides with $L_{\cM}^{p'}(X_1;L^{q'}(X_2;H^*))$.

\subsection{Proof of the sufficiency in Lemma 2.1}

In order to complete the proof of Lemma 2.1, it remains to prove its
``if" part, which is the main concern in this subsection.

Let $G:\cM\to H^*$ be a $(\mu_1\times\mu_2)$-continuous vector
measure of bounded variation. We want to show that there exists a
$\wt g\in L^1_\cM(X_1;L^1(X_2;H^*))$ such that
\bel{G}G(E)=\int_E\wt gd\m_1d\m_2,\qq\forall E\in\cM.\ee
Firstly, we
show that if $E_0\in \cM$ has a positive
$(\mu_1\times\mu_2)$-measure, then $G$ has a Bochner integrable
Radon--Nikod\'{y}m derivative on an $\cM$-measurable set $B$
satisfying $B\subset E_0$ and $(\mu_1\times\mu_2)(B)>0$.

Denote by $|G|$ the variation of $G$, which is a scalar measure (see
\cite[Definition 4 and Proposition 9 of Chapter 1, pp.2--3]{DU}). It
is easy to see that $|G|$ is a $(\mu_1\times\mu_2)$-continuous
$\dbR^+$-valued measure. Applying the Radon--Nikod\'{y}m Theorem (to
$|G|$ and $\mu_1\times\mu_2$), one can find an $\cM$-measurable
subset $B$ of $E_0$ and a positive integer $k$ such that $|G|(A)\leq
k(\mu_1\times\mu_2)(A)$ for all $A\in \cM$ with $A\subset B$. Define
a linear functional $\ell$ on the subspace $\cS$ of simple functions
in $L_{\cM}^p(X_1,L^q(X_2,H))$ as follows:
$$\ell(f) = \sum_{i=1}^n G(E_i\cap B)(x_i),$$
where
$$f=\sum_{i=1}^n x_iI_{E_i},\qq x_i\in H,\q1\le i\le n,$$
with $\{E_i,\;1\le i\le n\}$ being a partition of $X_1\times X_2$.
It follows that
$$\ba{ll}
 \ds|\ell(f)| = \Big| \sum_{i=1}^n G(E_i\cap B)(x_i)  \Big| = \Big|
\sum_{i=1}^n \frac{G(E_i\cap B)}{(\mu_1\times\mu_2)(E_i\cap
B)}\Big((\mu_1\times\mu_2)(E_i\cap B)x_i\Big) \Big|\\
\ns \ds\qq\;\; \leq \sum_{i=1}^n k |(\mu_1\times\mu_2)(E_i\cap
B)x_i| \leq k\|f\|_{L^1(X_i\times X_2;H)} \leq k
\mu_1(X_1)^{\frac{1}{p}}\mu_2(X_2)^{\frac{1}{q}}
\|f\|_{L^p_{\cM}(X_1;L^q(X_2;H))}.
 \ea$$
Therefore $\ell$ is a bounded linear functional  on $\cS$. By the
Hahn-Banach Theorem, it has a bounded linear extension to
$L_{\cM}^p(X_1,L^q(X_2,H))$ (The extension is still denoted by
$\ell$). Hence there exists a $g\in
L_{\cM}^{p'}(X_1,L^{q'}(X_2,H^*))$ such that
$$\ell(f) = \int_{X_1\times X_2} \langle f,g \rangle_{H,H^*}d\mu_1
d\mu_2 \qq \forall\; f\in L_{\cM}^p(X_1,L^q(X_2,H)).$$
We have
$$G(E\cap B)(x)=\ell(xI_{E})=\int_E\langle x,g
\rangle_{H,H^*}d\mu_1 d\mu_2,\qq\forall x\in H,~E\in M.$$
Since $g\in L_{\cM}^{p'}(X_1,L^{q'}(X_2,H^*))$ is Bochner
integrable, we see that
$$G(E\cap B)(x)=\Big( \int_E g d\mu_1
d\mu_2 \Big)(x),\qq\forall x\in H,\q E\in M.$$
Consequently,
\bel{oki1}\ds G(E\cap B) = \int_E g d\mu_1 d\mu_2,\qq\forall\;E\in
\cM. \ee
Noting that $B\in \cM$, and therefore replacing $E$ in (\ref{oki1})
by $E\cap B$, we see that
$$G(E\cap B) = \int_{E\cap B} g d\mu_1d\mu_2,\qq\forall\;E\in \cM.$$

Now by the Exhaustion Lemma (\cite[page 70]{DU}), there exist a
sequence $\ds\{A_n\}_{n=1}^{\infty}$ of disjoint members of $\cM$
such that $\ds\bigcup_{n=1}^{\infty}A_n = X_1\times X_2$ and a
sequence $\ds\{g_n\}_{n=1}^{\infty}$ of Bochner integrable functions
on $X_1\times X_2$ such that
$$G(E\cap A_n) = \int_{E\cap A_n}
g_n d\mu_1 d\mu_2,\qq\forall E\in \cM,\q n\in \dbN.$$
Define $\widetilde{g}: X_1\times X_2 \to H^*$ by
$\widetilde{g}(x_1,x_2) = g_n(x_1,x_2)$ if $(x_1,x_2)\in A_n$. It is
obvious that $\widetilde{g}$ is $(\mu_1\times \mu_2)$-measurable.
Moreover, for each $E\in \cM$ and all $m\in \dbN$, it holds
$$G\Big(E\bigcap\Big(\bigcup_{n=1}^m A_n\Big)\Big) = \int_E
\widetilde{g}I_{\cup_{n=1}^m A_n} d\mu_1 d\mu_2.$$
Consequently,
$$G(E) = \lim_{m\to\infty}\int_E
\widetilde{g}I_{\cup_{n=1}^m A_n} d\mu_1 d\mu_2,\qq\forall E\in
\cM.$$
For $h\in H^{**}$, the variation
$$|G(h)|(X_1\times X_2) \geq \lim_{m\to\infty}\int_{X_1\times X_2}
|\langle h , \widetilde{g}\rangle_{H^{**},H^*}|I_{\cup_{n=1}^m A_n}
d\mu_1 d\mu_2.$$
Hence by the Monotone Convergence Theorem, $\langle h,\widetilde{g}
\rangle_{H^{**},H^*}\in L^1_{\cM}(X_1;L^1(X_2;\dbR))$ for each $h\in
H^{**}$. If $E\in \cM$ and $h\in H^{**}$, from the Dominate
Convergence Theorem, we have
\begin{eqnarray*}
&\,& \langle h, G(E) \rangle_{H^{**},H^*} =
\lim_{m\to\infty}\int_{X_1\times X_2}  \langle h,\widetilde{g}
\rangle_{H^{**},H^*}I_{\cup_{n=1}^m A_n} d\mu_1 d\mu_2\\
&\,& \qq\qq\qq\q\; = \int_{X_1\times X_2}  \langle h,\widetilde{g}
\rangle_{H^{**},H^*}  d\mu_1 d\mu_2.
\end{eqnarray*}
Therefore $\widetilde{g}$ is Pettis integrable and its Pettis
integration P-$\ds\int_{X_1\times X_2} \widetilde{g} d\mu_1 d\mu_2 =
G(E)$ for each $E\in \cM$. Since $|G|(X_1\times X_2)$ is finite,
$\ds\int_{X_1\times X_2}|\widetilde{g}|_{H^*}I_{\cup_{n=1}^m A_n}
d\mu_1 d\mu_2 \leq |G|(X_1\times X_2)$ for all $m\in\dbN$. By the
Monotone Convergence Theorem, $|\widetilde{g}|_{H^*}\in
L^1_{\cM}(X_1;L^1(X_2;\dbR))$. Hence $\widetilde{g}$ is Bochner
integrable. Since the Pettis and Bochner integrals coincide whenever
they coexist, we obtain (\ref{G}), proving the Radon-Nikod\'ym
property of $H^*$ with respect to $(X_1\times
X_2,\cM,\m_1\times\m_2)$.

\subsection{A corollary of Lemma 2.1}

We now look an interesting corollary of Lemma 2.1. We first state
the following.

\ms

\bf Lemma 2.2. \sl Let
\bel{}\cM=\Big\{A\in\cB[0,T]\otimes\cF_T\bigm|t\mapsto
I_A(t,\cd)\hb{ is $\dbF$-progressively measurable }\Big\}.\ee
Then $\cM$ is a sub-$\si$-field of $\cB[0,T]\otimes\cF_T$. Moreover,
a process $\f:[0,T]\times\O\to H$ is $\dbF$-progressively measurable
if and only if it is $\cM$-measurable.

\ms

\bf Remark 2.3. \rm It is easy to see that the same conclusion in
Lemma 2.2 holds for any given filtration $\dbF$ (i.e., it is not
necessarily the natural filtration generated by the Brownian motion
$\{W(t)\}_{t\ge0}$), and also if one replaces the $\dbF$-progressive
measurability by any other measurability requirement, for examples,
adapted, optional or predictable, etc.

\ms

According to Lemmas 2.1 and 2.2, we have the following interesting
corollary, whose proof is straightforward.

\ms

\bf Corollary 2.4. \sl Let $0<s\leq T$ and $H^*$ have the
Radon--Nikod\'{y}m property with respect to
$([0,T]\times\O,\cM,m\times\dbP)$ (where $m$ is the Lebesgue
measure). Then the following identities hold:
\bel{}\left\{\ba{ll}
\ns\ds
L^p_\dbF(\O;L^q(0,s;H))^*=L^{p'}_\dbF(\O;L^{q'}(0,s;H^*)),\\
\ns\ds
L^q_\dbF(0,s;L^p(\O;H))^*=L^{q'}_\dbF(0,s;L^{p'}(\O;H^*)).\ea\right.\qq1\le
p,q<\infty.\ee

\rm

\ms

The above is a Riesz-type Representation Theorem for the dual of
spaces $L^p_\dbF(\O;L^q(0,s;H))$ and $L^q_\dbF(0,s;L^p(\O;H))$,
which will be very useful below.

\ms

We refer to \cite{L} for an application of Corollary 2.4 in the
study of null controllability of forward stochastic heat equations
with one control. We will give more application of this result in
our forthcoming papers;

\section{Answers to Problems (E) and (R)}\label{s3}

In this section, we return to our complete filtered probability
space $(\O,\cF,\dbF,\dbP)$ and give answers to Problems (E) and (R).

For any $p\in[1,\infty)$ and $0<s\leq T$, define an operator
$\dbL_s:L^p_\dbF(\O;L^1(0,s;H))\to L^p_{\cF_s}(\O;H)$  by
$$\dbL_s\big(u(\cd)\big)=\int_0^su(t)dt,\qq\forall\; u(\cd)\in
L^p_\dbF(\O;L^1(0,s;H)).$$
Concerning Problem (E), noting that
$L^1_\dbF(0,s;L^p(\O;H))\subseteq L^p_\dbF(\O;L^1(0,s;H))$, we give
the following positive answer (which is a little stronger than the
desired (\ref{1.11})):

\ms

\bf Theorem 3.1. \sl If $H^*$ has the Radon-Nikod\'ym property, then
\bel{2.1}\dbL_s\(L^1_\dbF(0,s;L^p(\O;H))\)=L^p_{\cF_s}(\O;H).\ee
Moreover, for each $\phi(\cd,s)\in L^p_{\cF_s}(\O;H)$, there is a
$\varsigma(\cd,s)\in L^1_\dbF(0,s;L^p(\O;H))$ such that
 \bel{kjld}
 \left\{
 \ba{ll}
 \ds
 \dbL_s\big(\varsigma(\cd,s)\big)=\phi(\cd,s),\\
 \ns
 \ds \|\varsigma(\cd,s)\|_{
L^1_\dbF(0,s;L^p(\O;H))}\le
\|\phi(\cd,s)\|_{L^1_\dbF(0,s;L^p(\O;H))}.
 \ea
 \right.
 \ee
(In general, the above $\varsigma(\cd,s)$ is NOT unique.)

\ms

\rm

The result in Theorem 3.1 turns out to be sharp for
$p\in(1,\infty)$. Indeed, we have the following result of negative
nature.

\ms

\bf Theorem 3.2. \sl For any $p\in(1,\infty)$ and any $r\in
(1,\infty]$, it holds that
\bel{oo2o1}\dbL_s\(L^r_\dbF(0,s;L^p(\O;H))\)\subsetneq
L^p_{\cF_s}(\O;H).\ee

\ms

\bf Remark 3.3. \rm 1) In \cite[VI, 68, pp. 130--131]{DM} and
\cite{D}, some Radon-Nikod\'ym type theorems were established for
real-valued or vector-valued processes with finite variation.
However, it seems that none of these results could be applied to
prove Theorem 3.1.

\ms

2) Thanks to Remark 2.3, the conclusion in Theorem 3.1 holds for any
given filtration $\dbF$; and one may replace the $\dbF$-progressive
measurability by any other measurability requirement.

\ms

3) We believe that (\ref{2.1}) is sharp in the sense that, for any
$r\in (1,\infty]$ and any $p\in [1,\infty]$,
\bel{oo2.1} \left\{ \ba{ll}
\dbL_s\(L^r_\dbF(0,s;L^p(\O;H))\)\subsetneq
L^p_{\cF_s}(\O;H),\\[2mm]
\dbL_s\(L^p_\dbF(\O;L^r(0,s;H))\)\subsetneq L^p_{\cF_s}(\O;H).
 \ea
 \right.\ee
Theorem 3.2 shows that the first conclusion in (\ref{oo2.1}) is true
for $p\in(1,\infty)$, and that, noting (\ref{1.3}), the second
conclusion in (\ref{oo2.1}) is true for $p\in(1,r]\cap (1,\infty)$.
The general case is under our investigation. Note that the above can
also be written as
\bel{}
\dbL_s\left(\bigcup_{q>1}L_\dbF^p(\O;L^q(0,s;H))\right)\subsetneq
L^p_{\cF_s}(\O;H).
\ee

As a consequence of Theorem 3.1 and the Martingale Representation
Theorem, our answer to Problem (R) is as follows:

\ms

\bf Corollary 3.4. \sl If $H$ is a Hilbert space, then for any
$p\in[1,\infty)$, one can find a constant $C>0$ such that for any
$\z(\cd)\in L^p_\dbF(\O;L^2(0,T;H))$, there is a $u(\cd)\in
L^1_\dbF(0,T;L^p(\O;H))$ so that equality $(\ref{R})$ holds and
 \bel{kj0ld}
 \|u(\cd)\|_{L^1_\dbF(0,T;L^p(\O;H))}
\le C\|\z(\cd)\|_{L^p_\dbF(\O;L^2(0,T;H))}.
  \ee

\ms

\bf Remark 3.5. \rm By point 2) in Remark 3.3, it is easy to see
that the conclusion in Corollary 3.4 holds also for adapted or
optional or predictable stochastic processes.

\ms

Corollary 3.4 shows the existence for the representation of It\^o
integrals by Lebesgue/Bochner integrals. The proof of Corollary 3.4
follows easily from Theorem 3.1 by noting the well-known result that
any Hilbert space has the Radon-Nikod\'{y}m property (e.g.,
\cite{DU}) and using also the Burkholder-Davis-Gundy inequality for
vector-valued stochastic processes (see \cite[Theorem 5.4]{CV} and
\cite[Corollary 3.11]{Neerven}). The rest of this section is devoted
to proving Theorems 3.1--3.2.

\ms

In order to prove Theorems 3.1--3.2, besides Corollary 2.4, we need
the following result concerning range inclusion for operators, which
can be found in \cite[Lemma 4.13, pp. 94--95 and Theorem 4.15, p.
97]{R}, for example.

\ms

\bf Lemma 3.6. \sl Suppose $B_X$ and $B_Z$ are the open unit balls
in Banach spaces $X$ and $Z$ , respectively. Let $L:\;X\to Z$ be a
linear bounded operator whose range is denoted by $\cR(L)$, and
whose adjoint operator is denoted by $L^*:Z^*\to X^*$. Then, the
following two conclusions hold

\ms

{\rm(i)} If $\cR (L)=Z$, then there is a constant $C>0$ such that
\bel{123}\|z^*\|_{Z^*}\le C\|L^*z^*\|_{X^*},\qq\forall\;z^*\in
Z^*.\ee

\ms

{\rm(ii)} If $(\ref{123})$ holds for some constant $C>0$, then
\bel{123-1}B_Z\subset CL(B_X)\equiv\big\{CLx\;\big|\;x\in
B_X\big\}.\ee

\ms

\bf Remark 3.7. \rm 1) Clearly, by Lemma 3.6, we see that $\cR
(L)=Z$ if and only if (\ref{123}) holds for some constant $C>0$. But
this lemma goes a little further than this. Indeed, the second
conclusion of this lemma provides a ``quantitative" characterization
$B_Z\subset CL(B_X)$, which is more delicate than $\cR (L)=Z$. We
shall use this result essentially when we answer Problem (C) in the
next section;

2) One should compare Lemma 3.6 with the following general range
inclusion result (e.g., \cite[Lemma 2.4 in Chap. 7]{LY}): Let $X,Y$
and $Z$ be Banach spaces with $X$ being reflexive, and both $F:Y\to
Z$ and $G:X\to Z$ be linear bounded operators. Then,
 \bel{olko}
 \ba{ll}
 |F^*z^*|_{Y^*}\le C|G^*z^*|_{X^*},\q\forall z^*\in Z^*,\hb{ \rm for some
constant }C>0\\
\iff \cR(F)\subseteq\cR(G).
 \ea
 \ee
As shown in \cite{B}, the equivalence (\ref{olko}) may fail whenever
$X$ is not reflexive. Nevertheless, when $F$ is surjective (in
particular when $Y=Z$ and $F = I$, the identity operator, the case
considered in Lemma 3.6), this equivalence remains to be true (even
without the reflexivity assumption for $X$) (see \cite[Theorem 1.2
and Remark 1.3]{A}). We refer to \cite{WZ} for further range
inclusion results.

\ms

Further, we need the following property for Wiener integrals,  a
special case of It\^o integrals with deterministic integrands (e.g.,
\cite[Theorem 2.3.4 in Chapter 2, p. 11]{Kuo}).

\ms

\bf Lemma 3.8. \sl For each $0\le a<b\le T$ and $f \in L^2(a, b)$
(for which $f$ is a deterministic function, i.e., it does not depend
on $\o\in\O$), the Wiener integral $\int^b_a f(t) dW(t)$ is a
Gaussian random variable with mean $0$ and variance $\int^b_a\left|
f(t)\right|^2 dt$.

\ms

\rm

We are now in a position to prove Theorems 3.1--3.2.

\ms

\it Proof of Theorem 3.1. \rm It suffices to show (\ref{kjld}).
Since $L^1_\dbF(0,s;L^p(\O;H))\subseteq L^p_\dbF(\O;L^1(0,s;H))$
(algebraically and topologically), the restriction of operator
$\dbL_s:L^p_\dbF(\O;L^1(0,s;H))\to L^p_{\cF_s}(\O;H)$ to
$L^1_\dbF(0,s;L^p(\O;H))$ is a bounded linear operator from
$L^1_\dbF(0,s;L^p(\O;H))$ to $L^p_{\cF_s}(\O;H)$ (For simplicity, we
still denote it by $\dbL_s$). By Conclusion (ii) in Lemma 3.6 and
Corollary 2.4, by a simple scaling, we see that the desired result
(\ref{kjld}) is implied by the following:
\bel{2.2}\|\dbL_s^*\eta\|_{L^\infty_\dbF(0,s;L^{p'}(\O;H^*))}\ge\|\eta\|_{L^{p'}_{\cF_s}(\O;H^*)},\qq\forall\;
\eta\in L^{p'}_{\cF_s}(\O;H^*).\ee

In order to prove (\ref{2.2}), let us first find the adjoint
operator $\dbL_s^*$ of $\dbL_s$. For any $u(\cd)\in
L^1_\dbF(0,s;L^p(\O;H))$, and $\eta\in
L^p_{\cF_s}(\O;H)^*=L^{p'}_{\cF_s}(\O;H^*)$, we have
\bel{}\ba{ll}
\ns\ds\lan\dbL_s
u,\eta\ran=\dbE\left(\int_0^s u(t)dt,\eta\right)_{H,H^*}=\int_0^s\dbE \Big(u(t),\eta\Big)_{H,H^*} dt\\
\ns\ds\qq\qq=\int_0^s\dbE\(u(t),\dbE[\eta\,|\,\cF_t]\)_{H,H^*}dt=\lan
u,\dbL_s^*\eta\ran,\ea\ee
which leads to
\bel{L*}\left\{\ba{ll}
\ns\ds\dbL_s^*:L^{p'}_{\cF_s}(\O;H^*)\to
L^1_\dbF(0,s;L^p(\O;H))^*=L^\infty_\dbF(0,s;L^{p'}(\O;H^*)),\\
\ns\ds(\dbL_s^*\eta)(t)=\dbE[\eta\,|\,\cF_t],\qq
t\in[0,s],~\forall\;\eta\in L^{p'}_{\cF_s}(\O;H^*).\ea\right.\ee
This gives a representation of the adjoint operator $\dbL_s^*$ of
$\dbL_s$.

\ms

Now, we let $p>1$. Making use of (\ref{L*}), we find that
\bel{2.15}\ba{ll}
\ns\ds\|\dbL_s^*\eta\|_{L^\infty_\dbF(0,s;L^{p'}(\O;H^*))}
=\left[\sup_{t\in[0,s]}\dbE\Big|\dbE[\eta\,|\,\cF_t]\Big|_{H^*}^{p'}\right]^{1\over
p'}\\
\ns\ds\ge\left[\dbE\Big|\dbE[\eta\,|\,\cF_s]\Big|_{H^*}^{p'}\right]^{1\over
p'}=\left[\dbE|\eta|^{p'}\right]^{1\over
p'}=\|\eta\|_{L^{p'}_{\cF_s}(\O;H^*)}.\ea\ee
Therefore, (\ref{2.2}) holds for $p>1$.

\ms

Next, for $p=1$, we have that
\bel{2.16}\ba{ll}
\ns\ds\|\dbL_s^*\eta\|_{L^\infty_\dbF(\O;L^\infty(0,s;H^*))}
=\esssup_{\o\in\O}\left[\sup_{t\in[0,s]}|\dbE[\eta\,|\,\cF_t]|_{H^*}\right]\\
\ns\ds\qq\qq\qq\qq\qq\ge\esssup_{\o\in\O}\[|\dbE[\eta\,|\,\cF_s]|_{H^*}\]=\esssup_{\o\in\O}|\eta(\o)|_{H^*}=\|\eta\|_{L^\infty_{\cF_s}(\O;H^*)}.\ea\ee
This implies that our conclusion also holds for $p=1$.\endpf

\ms

\it Proof of Theorem 3.2. \rm Noting (\ref{basd}), it suffices to
prove Theorem 3.2 for $r\in (1,\infty)$. We use the contradiction
argument. Assume that
 \bel{kkqq}
 \dbL_s\(L^r_\dbF(0,s;L^p(\O;H))\)=
L^p_{\cF_s}(\O;H),\q \hb{ for some }p, r\in (1,\infty).
 \ee
Since $L^r_\dbF(0,s;L^p(\O;H))\subseteq
L^1_\dbF(0,s;L^p(\O;H))\subseteq L^p_\dbF(\O;L^1(0,s;H))$
(algebraically and topologically), the restriction of operator
$\dbL_s:L^p_\dbF(\O;L^1(0,s;H))\to L^p_{\cF_s}(\O;H)$ to
$L^r_\dbF(0,s;L^p(\O;H))$ is again a bounded linear operator from
$L^r_\dbF(0,s;L^p(\O;H))$ to $L^p_{\cF_s}(\O;H)$ (For simplicity, we
still denote it by $\dbL_s$). Similar to (\ref{L*}), the
representation of the adjoint operator $\dbL_s^*$ of $\dbL_s$ is
given as follows:
\bel{L*1}\left\{\ba{ll}
\ns\ds\dbL_s^*:\ L^{p'}_{\cF_s}(\O;H^*)\to
L^{r'}_\dbF(0,s;L^{p'}(\O;H^*)),\\
\ns\ds(\dbL_s^*\eta)(t)=\dbE[\eta\,|\,\cF_t],\qq
t\in[0,s],~\forall\;\eta\in L^{p'}_{\cF_s}(\O;H^*).\ea\right.\ee

By (\ref{kkqq}), using the first conclusion in Lemma 3.6 and noting
Corollary 2.4, we conclude that there exists a constant $C>0$ such
that for any $\eta\in L^{p'}_{\cF_s}(\O;H^*)$, it holds that
\begin{equation}\label{range inequality}
\|\eta\|_{L^{p'}_{\cF_s}(\O;H^*)}\leq C\|\dbL_s^*\eta\|_{
L^{r'}_\dbF(0,s;L^{p'}(\O;H^*))},
\end{equation}
where $r'=r/( r-1)$.

Fix any $x_0\in H^*$ satisfying $|x_0|_{H^*} =1$ (which is
independent of the time variable $t$ and the sample point $\omega$).
Consider a sequence of random variables $\{\eta_n\}_{n=1}^\infty$
defined by
$$
\eta_n = \int_0^s e^{nt}dW(t)x_0,\qq n\in \mathbb{N}.
$$
It is obvious that $\eta_n \in L^{p'}_{\cF_s}(\O;H^*)$ for any $n\in
\mathbb{N}$. By Lemma 3.8, the integral $\displaystyle\int_0^s
e^{nt}dW(t)$ is a Gaussian random variable with mean $0$ and
variance $\frac{e^{2ns}-1}{2n}$. Hence,
\begin{equation}\label{ra0}
 \ba{ll}\ds
\ds\left[\mathbb{E}\left|\int_0^s
e^{nt}dW(t)\right|^{p'}\right]^{\frac{1}{p'}}\3n&\ds =
\left[\int_{-\infty}^{\infty}\frac{\sqrt{n}\left|x\right|^{p'}}{\sqrt{(e^{2ns}-1)\pi}}e^{-\frac{nx^2}{e^{2ns}-1}}dx\right]^{\frac{1}{p'}}\\[3mm]
 &\ds=
\left[\int_{-\infty}^{\infty}\left(\frac{e^{2ns}-1}{n}\right)^{p'/2}\frac{\left|x\right|^{p'}}{\sqrt{\pi}}e^{-x^2}dx\right]^{\frac{1}{p'}}
\\[3mm]
 &\ds=\left(\frac{1}{\sqrt{\pi}}\int_{-\infty}^{\infty}\left|x\right|^{p'}e^{-x^2}dx\right)^{\frac{1}{p'}}\sqrt{\frac{e^{2ns}-1}{n}}.
 \ea
\end{equation}

Now, by (\ref{ra0}), it is easy to see that
\begin{equation}\label{ra1}
\ba{ll}\ds
\|\eta_n\|_{L^{p'}_{\cF_s}(\O;H^*)}\3n&\ds=\left[\mathbb{E}\left|\int_0^s
e^{nt}dW(t)x_0\right|^{p'}\right]^{\frac{1}{p'}}=\left[\mathbb{E}\left|\int_0^s
e^{nt}dW(t)\right|^{p'}\right]^{\frac{1}{p'}} \\[3mm]
 &\ds=
\left(\frac{1}{\sqrt{\pi}}\int_{-\infty}^{\infty}\left|x\right|^{p'}e^{-x^2}dx\right)^{\frac{1}{p'}}\sqrt{\frac{e^{2ns}-1}{n}}.
 \ea
\end{equation}
Using (\ref{ra0}) again, we have
 \bel{ra2}
 \ba{ll}
\ds
\big\|\mathbb{E}[\eta_n|\cF_t]\big\|_{L^{r'}_{\dbF}(0,s;L^{p'}(\O;H^*))}
\3n&\ds= \left\{\int_0^s \left[\mathbb{E}\left|\int_0^t
e^{n\tau}dW(\tau)x_0\right|^{p'}\right]^{\frac{r'}{p'}}dt\right\}^{\frac{1}{r'}}\\
 \ns
&\ds= \left\{\int_0^s \left[\mathbb{E}\left|\int_0^t
e^{n\tau}dW(\tau)\right|^{p'}\right]^{\frac{r'}{p'}}dt\right\}^{\frac{1}{r'}}\\
 \ns
&\ds =\left\{\int_0^s \left[
\left(\frac{1}{\sqrt{\pi}}\int_{-\infty}^{\infty}\left|x\right|^{p'}e^{-x^2}dx\right)^{\frac{1}{p'}}\sqrt{\frac{e^{2nt}-1}{n}}
\right]^{r'}dt\right\}^{\frac{1}{r'}}\\
 \ns
&\ds \leq\frac{1}{\sqrt{n}}
\left(\frac{1}{\sqrt{\pi}}\int_{-\infty}^{\infty}\left|x\right|^{p'}e^{-x^2}dx\right)^{\frac{1}{p'}}\left(\int_0^s
e^{nr't}dt\right)^{\frac{1}{r'}}\\
 \ns
&\ds \le\frac{1}{\sqrt{n}}
\left(\frac{1}{\sqrt{\pi}}\int_{-\infty}^{\infty}\left|x\right|^{p'}e^{-x^2}dx\right)^{\frac{1}{p'}}\frac{e^{ns}}{\left(nr'\right)^{\frac{1}{r'}}}.
\ea \ee
From (\ref{ra1}) and (\ref{ra2}), it follows that
$$
\lim_{n\to\infty}\frac{\big\|\mathbb{E}[\eta_n|\cF_t]\big\|_{L^{r'}_{\dbF}(0,s;L^{p'}(\O;H^*))}}
{\|\eta_n\|_{L^{p'}_{\cF_s}(\O;H^*)}}\le
\lim_{n\to\infty}\frac{e^{ns}}{\left(nr'\right)^{\frac{1}{r'}}\sqrt{e^{2ns}-1}}=0.
 $$
This, combined with (\ref{L*1}), gives
$$\lim_{n\to\infty}\frac{\big\|\dbL_s^*\eta_n\big\|_{L^{r'}_{\dbF}
(0,s;L^{p'}(\O;H^*))}}{\|\eta_n\|_{L^{p'}_{\cF_s}(\O;H^*)}}=0,$$
which contradicts inequality (\ref{range inequality}). This
completes the proof of Theorem 3.2. \endpf

\section{Answer to Problem (C)}\label{s4}

This section is addressed to give a positive answer to Problem (C).

Theorem 3.1 tells us that any It\^o integral $\ds\int_0^s\z(t)dW(t)$
with $\z(\cd)\in L^p_\dbF(\O;L^2(0,T;H))$ admits a (parameterized)
Bochner integral representation, i.e. we can find a representor
$u(\cd,s)\in L^1_\dbF(0,s;L^p(\O;H))$ (which is of course NOT
unique) such that
\begin{equation}\label{repid1}
\int_0^s\z(t)dW(t) = \int_0^s u(t,s)dt,\qq \forall\; s\in [0,T].
\end{equation}
Put $Z\equiv L^1_\dbF(0,T;L^p(\O;H))$. We now show that one can
choose a $u(\cd,s)$, which is continuous in $Z$ with respect to $s$,
such that (\ref{repid1}) holds. More precisely, we have the
following result:

\ms

\bf Theorem 4.1. \sl For any given $\z(\cd)\in
L_{\mathbb{F}}^p(\O;L^2(0,T;H))$, define a (set-valued) mapping $F:
[0,T]\to 2^Z$ by
\begin{eqnarray}\label{F}
F(s)=\Big\{  \eta(\cd,s)\in Z\;\Big| \int_0^s \eta(t,s)dt = \int_0^s
\z(t)dW(t), \hb{ and }\eta(t,s)=0, \,\forall\, t>s\Big\},\;\;
\forall\; s\in [0,T].
\end{eqnarray}
Then $F$ has a continuous selection $f$.

\ms

\bf Remark 4.2. \rm If we choose $u(\cd,s)$ to be the above $f(s)$,
then $u(\cd,s)$ is the desired process (for (\ref{repid1})), which
is continuous in $Z$ with respect to $s$.

\ms

Before proving Theorem 4.1, we recall the following useful
preliminary results.

\ms

\bf Lemma 4.3. \sl Let $X$ and $Y$ be two topological spaces. Then,
for any set-valued mapping $\phi:X\to 2^Y$,  the following two
statements are equivalent:

\ms

{\rm (i)} The map $\phi$ is lower semi-continuous, i.e., for any
open subset $V$ of $Y$, the set $\Big\{ x\in X\;\Big|\; \phi(x)\cap
V \neq \emptyset \Big\}$ is open in $X$;

\ms

{\rm(ii)} If $x\in X$, $y\in \phi(x)$, and $V$ is a neighborhood of
$y$ in $Y$, then there exists a neighborhood $U$ of $x$ in $X$ such
that for every $x'\in U$, there exists a $y'\in \phi(x')\cap V$.

\ms

\bf Lemma 4.4. {\rm (\cite[Theorem 3.2$\,''$]{Micheal})} \sl The
following properties of a $T_1$-space are equivalent:

\ms

{\rm(i)} $X$ is paracompact (i.e., any open cover of $X$ admits a
locally finite open refinement, which is the case if $X$ is compact
or is a metric space);

{\rm(ii)} If $Y$ is a Banach space, then every lower semi-continuous
mapping  $F:X\to  2^Y$ such that $F(x)$  is  a non-empty, closed,
convex subset of $Y$ for any $x\in X$, admits a continuous
selection, i.e., there exists a continuous mapping $f:X\to Y$ such
that $f(x)\in F(x)$ for any $x\in X$.

\ms

\rm

We can now give a proof of Theorem 4.1.

\ms

{\it Proof of Theorem 4.1}. The main idea is to use Lemma 4.4. It is
obviously that $[0,T]$ is an $T_1$-space and is paracompact. Hence
we need only to prove that $F(s)$ is a non-empty, closed, convex
subset of $Z$ for any $s\in [0,T]$ and $F$ is lower semi-continuous.
By Theorem 3.1, we see that $F(s)$ is non-empty. Also, it is very
easy to check that $F(s)$ is a convex subset of $Z$ and is closed in
$Z$.

It remains to show that $F$ is lower semi-continuous. Fix any $s\in
[0,T]$, any $\eta(\cdot,s)\in F(s)$, and any neighborhood $V$ of
$\eta(\cdot,s)$ in $Z$. Clearly, there exists a $\d>0$ such that
$$
V_1 = \big\{ z(\cd)\in Z\,\big|\,\|z(\cd)-\eta(\cdot,s)\|_Z < \d
\big\}\subset V.
$$
We claim that there exists an $\varepsilon>0$ such that for any $r$
satisfying $|r-s|<\varepsilon$, it holds that
 \bel{ok1}
 F(r)\cap V_1
\neq\emptyset.
 \ee
This claim will yield the lower semi-continuity of $F(\cd)$. To
prove out claim, we first make use of the Burkholder-Davis-Gundy
inequality for vector-valued stochastic process (see \cite[Theorem
5.4]{CV} and \cite[Corollary 3.11]{Neerven}) to get the following:
\begin{equation}\label{BDG}
\dbE\Big|\int_r^s \z(t)dW(t)\Big|_{H}^p \leq  \dbE \Big[ \sup_{r\leq
h\leq s }\Big|\int_r^h \z(t)dW(t)\Big|_H^p \Big] \leq
C\dbE\Big[\int_r^s |\z(t)|_{H}^2dt\Big]^{\frac{p}{2}}.
\end{equation}
Choose an increasing sequence $\{r_k\}_{k=1}^\infty$ such that $0\le
r_1\le r_2\le \cdots \le r_k\le r_{k+1}\le \cdots\to s$. Since
$\z(\cd)\in L^p_\dbF(\O;L^2(0,T;H))$, by the Dominated Convergence
Theorem, we have
$$\lim_{k\to \infty}\dbE\Big[\int_{r_k}^s
|\z(t)|_{H}^2dt\Big]^{\frac{p}{2}}=\lim_{{k\to
\infty}}\dbE\Big[\int_0^T\chi_{[{r_k},s]}
|\z(t)|_{H}^2dt\Big]^{\frac{p}{2}}=0.$$
Hence,
\bel{yu1}\lim_{r\to s}\dbE\Big[\int_r^s
|\z(t)|_{H}^2dt\Big]^{\frac{p}{2}}\le \lim_{k\to
\infty}\dbE\Big[\int_{r_k}^s |\z(t)|_{H}^2dt\Big]^{\frac{p}{2}}=0.
\ee
Therefore, it follows from (\ref{BDG}) that there exists an
$\varepsilon_1>0$ such that for any $0\leq s-r<\varepsilon_1$, the
following holds
\bel{the1}\Big\|\int_r^s \z(t)dW(t)\Big\|_{L_{\cF_s}^p(\O;H)}
<\frac{\d}{3}.
 \ee
On the other hand, by the H\"older inequality and using the
Dominated Convergence Theorem, similar to the proof of (\ref{yu1}),
we see that there exists an $\varepsilon_2>0$ (may depend on $s$)
such that for any $0\leq s-r<\varepsilon_2$, it holds
 \bel{the2}
 \Big\|\int_r^s \eta(t,s)dt\Big\|_{L_{\cF_s}^p(\O;H)}\le\int_r^s
 \|\eta(t,s)\|_{L^p_{\cF_s}(\O;H)}dt=\int_r^s
\Big[\dbE\big|\eta(t,s)\big|^{p}_{H}\Big]^{\frac{1}{p}}dt<\frac{\d}{3}.
 \ee
Put $\varepsilon_3 = \min\{\varepsilon_1,\varepsilon_2\}$. From
(\ref{the1})--(\ref{the2}) and noting that $\ds\int_0^s \eta(t,s)dt
= \int_0^s \z(t)dW(t)$, we conclude that for any $r$ satisfies
$0\leq s-r<\varepsilon_3$, it holds that
 \bel{ine1}
 \ba{ll}
\ds\Big\|\int_0^r \eta(t,s)dt-\int_0^r
\z(t)dW(t)\Big\|_{L_{\cF_s}^p(\O;H)}\\
\ns \ds\le\2n\Big\|\int_0^r \eta(t,s)dt-\int_0^s
\eta(t,s)dt\Big\|_{L_{\cF_s}^p(\O;H)}\1n +\1n \Big\|\int_0^r
\z(t)dW(t)-\int_0^s
\z(t)dW(t)\Big\|_{L_{\cF_s}^p(\O;H)}<\frac{2\d}{3}.
 \ea
 \ee
By the second conclusion in Theorem 3.1 and noting (\ref{ine1}), we
see that there is a $\ds\phi(\cdot,r)\in L^1_\dbF(0,r;L^p(\O;H))$
such that
$\ds\|\phi(\cdot,r)\|_{L^1_\dbF(0,r;L^p(\O;H))}<\frac{2\d}{3}$, and
$$
\int_0^r \phi(t,r)dt =\int_0^r \z(t)dW(t) - \int_0^r \eta(t,s)dt.
$$
Put $\varrho(\cdot,r) = \chi_{[0,r]}\phi(\cdot,r) +
\chi_{[0,r]}\eta(\cdot,s)$. It is obvious that $\varrho(\cdot,r)\in
F(r)$, and
$$\Big\| \eta(\cdot,s) -\varrho(\cdot,r)
\Big\|_{L_{\dbF}^1(0,s;L^p(\O,H))} \leq \int_r^s
\Big[\dbE|\eta(t,s)|^{p}_{H}\Big]^{\frac{1}{p}}dt+\|\phi(\cdot,r)\|_{L^1_\dbF(0,r;L^p(\O;H))}<\d.$$
Therefore, for any $0\leq s-r <\varepsilon_3$, it holds that
$\varrho(\cdot,r)\in V_1$, which gives (\ref{ok1}). By a similar
argument, one can show that there exists an $\varepsilon_4>0$ such
that (\ref{ok1}) holds for any $0\leq r-s <\varepsilon_4$. Choosing
$\varepsilon = \min\{ \varepsilon_3,\varepsilon_4 \}$, we see that
(\ref{ok1}) holds for any $|r-s|<\varepsilon$. By Lemma 4.3, we know
that $F:[0,T] \to Z$ is lower semi-continuous.

Finally, thanks to Lemma 4.4, we conclude that there exists a
continuous selection $f$ of $F$. \endpf

\section{Two Illustrative Applications}\label{s5}

In this section, we give two simple applications of our Theorems
3.1--3.2. More interesting and sophisticated applications will be
presented in our forthcoming publications.

\subsection{Application to the controllability problem}\label{5sub1}

Consider a one-dimensional controlled stochastic differential
equation:
\bel{1.23}dx(t)=[bx(t)+u(t)]dt+\si dW(t),\ee
with $b$ and $\si$ being given constants. We say that system
(\ref{1.23}) is {\it exactly controllable} if for any $x_0\in\dbR$
and $x_T\in L^p_{\cF_T}(\O;\dbR)$, there exists a control $u(\cd)\in
L^p_\dbF(\O;L^1(0,T;\dbR))$ such that the corresponding solution
$x(\cd)$ satisfies $x(0)=x_0$ and $x(T)=x_T$.
By variation of constant formula, we have
$$x(T)=e^{bT}x_0+\int_0^Te^{b(T-t)}u(t)dt+\int_0^Te^{b(T-t)}\si
dW(t).$$
Thus, exact controllability is equivalent to the following:
\bel{1.26}x_T-e^{bT}x_0-\int_0^Te^{b(T-t)}\si
dW(t)=\int_0^Te^{b(T-t)}u(t)dt.\ee
Since $x_T\in L^p_{\cF_T}(\O;\dbR)$, there exists a unique
$\z(\cd)\in L^p_\dbF(\O;L^2(0,T;\dbR))$ such that
$$x_T=\dbE x_T+\int_0^T\z(t)dW(t).$$
Hence, to ensure (\ref{1.26}), it suffices to have
$$\dbE x_T-e^{bT}x_0+\int_0^T\big[\z(t)-e^{b(T-t)}\si\big]
dW(t)=\int_0^Te^{b(T-t)}u(t)dt,$$
which is guaranteed by Theorem 3.1. This means that (\ref{1.23}) is
exactly controllable.

On the other hand, surprisingly, in virtue of \cite[Theorem
3.1]{Peng}, it is clear that system (\ref{1.23}) is NOT exactly
controllable if one restricts to use admissible controls $u(\cd)$ in
$L^2_\dbF(\O;L^2(0,T;\dbR))$! Moreover, by Theorem 3.2, we see that
system (\ref{1.23}) is NOT exactly controllable, either provided
that one uses admissible controls $u(\cd)$ in
$L^2_\dbF(\O;L^q(0,T;\dbR))$ for any $q\in(1,\infty]$. This leads to
a corrected formulation for the exact controllability of stochastic
differential equations, as presented below.

We consider the following linear stochastic differential equation:
 \bel{sd1}
 \left\{
 \ba{ll}
 dy(t) =\big[Ay(t)+Bu(t)\big]dt+\big[Cy(t)+Du(t)\big]dW(t), \qq 0 \le  t \le T, \\
 \ns
 y(0) = y_0\in \dbR^n,
 \ea
 \right.
 \ee
where $A,C \in \dbR^{n\times n}$ and $B,D \in \dbR^{n\times m}$
($n,m\in\dbN$) are matrices. Various controllability issues for
system (\ref{sd1}) were studied, say, in \cite{BQT,CLPY,G, Peng} and
the references cited therein. Note however that, unlike the
classical deterministic case, as far as we know, there exist no
universally accepted notions for controllability in the stochastic
setting so far.

Motivated by the above observation, we introduce the following:

\ms

\bf Definition 5.1. \rm System (\ref{sd1}) is said to be exactly
controllable if for any $y_0\in \dbR^n$ and  $y_T\in
L^p_{\cF_T}(\O;\dbR^n)$, there exists a control $u(\cd)\in
L^p_\dbF(\O;L^1(0,T;\dbR^m))$ such that $Du(\cd)\in
L^p_\dbF(\O;L^2(0,T;\dbR^n))$ and the corresponding solution
$y(\cd)$ of (\ref{sd1}) satisfies $y(T)=y_T$.

\ms

We need $Du(\cd)\in L^p_\dbF(\O;L^2(0,T;\dbR^n))$ in the above
definition because it appears in the It\^o integral
$\ds\int_0^T\big[Cy(t)+Du(t)\big] dW(t)$. It is clear that, for the
controllability of deterministic linear (time-invariant) ordinal
differential equations, there is no difference between the
controllability by using $L^1$ (in time) control and that by using
$L^2$ (or even analytic in time) control. However, our analysis
above indicates that things are completely different in the
stochastic setting. A detailed study of the controllability for
system (\ref{sd1}) (in the sense of Definition 5.1) seems to deviate
the theme of this paper, and therefore we shall address this topic
in our forthcoming works.

\subsection{Application to a Black-Scholes model}

 Consider a Black-Scholes market model
\bel{}\left\{\ba{ll}
\ns\ds dX_0(t)=rX_0(t)dt,\\
\ns\ds dX(t)(t)=bX(t)dt+\si X(t)dW(t),\ea\right.\ee
with $r,b,\si$ being constants. Under conditions of self-financing,
and no transaction costs, the investor's wealth process $Y(\cd)$
satisfies the following equation:
\bel{}dY(t)=\big[rY(t)+(b-r)Z(t)\big]dt+\si Z(t)dW(t),\ee
where $Z(t)$ is the amount invested in the stock. For convenience, a
European contingent claim with payoff at the maturity $T$ being
$\xi\in L^p_{\cF_T}(\O;\dbR)$ is identified with $\xi$. Any such a
$\xi$ is said to be {\it replicatable} if there exists a trading
strategy $Z(\cd)$ such that for some $Y_0$ (the price of the
contingent claim at $t=0$), one has
$$Y(0)=Y_0,\qq Y(T)=\xi.$$
In another word, contingent claim $\xi$ is replicatable if and only
if the following backward stochastic differential equation (BSDE,
for short) admits an adapted solution $(Y(\cd),Z(\cd))$:
\bel{}\left\{\ba{ll}
\ns\ds dY(t)=\big[rY(t)+(b-r)Z(t)\big]dt+\si Z(t)dW(t),\q
t\in[0,T],\\
\ns\ds Y(T)=\xi.\ea\right.\ee
In this case, $Y(t)$ is a price of the contingent claim at time $t$.
See \cite{EPQ} and \cite{Yong 2006} for some relevant presentations.
Now, let us look at an extreme case,
\bel{si=0}b-r>0,\qq\si=0.\ee
In this case, $\xi$ is replicatable if and only if the following
BSDE admits an adapted solution $(Y(\cd),Z(\cd))$:
$$\left\{\ba{ll}
\ns\ds dY(t)=\big[rY(t)+(b-r)Z(t)\big]dt,\qq t\in[0,T],\\
\ns\ds Y(T)=\xi.\ea\right.$$
Similar to the above subsection, we see that the above admits an
adapted solution $(Y(\cd),Z(\cd))$, which means that $\xi$ is
replicatable. Further, since $\xi$ is arbitrary, this also means
that the market with conditions (\ref{si=0}) is complete! This is a
little surprising since $\si=0$ in the market model. Some further
careful study along this line will be carried out in our future
investigations.

\end{document}